\documentclass[12pt]{article}

\usepackage{amssymb,amsmath,amsthm}

\newcommand{\vp}{\varepsilon}

\newcommand{\bb}[1]{{\mathbb{#1}}}

\newcommand{\N}{N}

\numberwithin{equation}{section}

\theoremstyle{plain}
\newtheorem{lem}{Lemma}[section]

\newtheorem{thm}[lem]{Theorem}
\newtheorem{cor}[lem]{Corollary}

\theoremstyle{definition}

\theoremstyle{remark}

\begin{document}

\begin{center}
\LARGE 
STRONG SINGULARITY OF SINGULAR MASAS\\
 IN  ${\rm{II}}_1$ FACTORS
\end{center}\bigskip

$$\begin{tabular}{cccc}

Allan M.~Sinclair&&& \hspace{.5in}Roger R.~Smith$^{*}$\\
\vspace{15pt}
Stuart A.~White &&& \hspace{.5in}Alan Wiggins\\
\vspace{-5pt}
School of Mathematics&&& \hspace{.5in}Department of Mathematics\\
\vspace{-5pt}
University of Edinburgh&&& \hspace{.5in}Texas A\&M University\\
\vspace{-5pt}
Edinburgh EH9 3JZ&&& \hspace{.5in}College Station, TX \ 77843\\
Scotland&&&\hspace{.5in}USA\\
\vspace{-5pt}
{\tt {A.Sinclair@ed.ac.uk}}&&& \hspace{.5in}{\tt {rsmith@math.tamu.edu}}\\

{\tt {s.a.white-1@ed.ac.uk}}&&& \hspace{.5in}{\tt {awiggins@math.tamu.edu}}\\

\end{tabular}$$
\vspace{.5in}

\abstract{
 A singular masa $A$ in a ${\rm{II}}_1$ factor $N$ is defined by the
property that any unitary $w\in N$ for which $A=wAw^*$ must lie in
$A$. A strongly singular masa $A$ is one that satisfies the
inequality
\[\|\bb E_A-\bb E_{wAw^*}\|_{\infty,2}\geq\|w-\bb E_A(w)\|_2\]
for all unitaries $w\in N$, where $\bb E_A$ is the conditional expectation of
$N$ onto $A$, and $\|\cdot\|_{\infty,2}$ is defined for bounded maps
$\phi :N\to N$ by $\sup\{\|\phi(x)\|_2:x\in N,\ \|x\|\leq 1\}$. Strong
singularity easily implies singularity, and the main result of this paper shows
the reverse implication.}

\vfill

$\underline{\hspace{1in}}$

\noindent $^*$Partially supported by a grant from the National Science
Foundation.
\newpage

\section{Introduction}\label{sec1}

\indent
In \cite{SS1}, the first two authors introduced the concept of a strongly
singular maximal abelian self--adjoint subalgebra (masa) in a 
${\rm{II}}_1$
factor $N$. For a bounded map $\phi:M\to N$ between any two finite von~Neumann
algebras with specified traces, the $\|\cdot\|_{\infty,2}$--norm is defined by
\begin{equation}\label{eq1.1}
\|\phi\|_{\infty,2}=\sup \{\|\phi(x)\|_2:x\in M,\ \|x\|\leq 1\}.
\end{equation}
A masa $A\subseteq N$ is then said to be {\emph{strongly singular}} when the
inequality
\begin{equation}\label{eq1.2}
\|\bb E_A-\bb E_{wAw^*}\|_{\infty,2}\geq \|w-\bb E_A(w)\|_2
\end{equation}
holds for all unitaries $w\in N$, where the notation $\bb E_B$ indicates the
unique trace preserving conditional expectation onto a von~Neumann subalgebra
$B$. Any
unitary which normalizes $A$ is forced,
by this relation, to lie in $A$, and so $A$ is singular, as defined in
\cite{Dix}.
The original purpose for introducing strong singularity was to have a metric
condition which would imply singularity, and which would be easy to verify in
a wide range of cases (see \cite{SS1,SS2} and the work of the third author,
\cite{W1}, on Tauer masas in the
hyperfinite $\rm{II}_1$ factor $R$).
Subsequently, \cite{PSS}, it was shown that every singular masa in
a separable  ${\rm{II}}_1$ factor (where this terminology indicates
norm--separability of the predual $N_*$) satisfies the following
weaker inequality, analogous to \eqref{eq1.2}:
\begin{equation}\label{eq1.3}
90\|\bb E_A-\bb E_{wAw^*}\|_{\infty,2}\geq \|w-\bb E_A(w)\|_2
\end{equation}
holds for all unitaries $w\in N$. This clearly suggested that every singular
masa should be strongly singular, and our objective in this paper is to prove
this result.


In 1983 Sorin Popa   introduced the
$\delta$--invariant for a masa in a ${\rm{II}}_1$ factor, \cite{Po_conv}. This
was
the first attempt to define a metric based invariant for a masa,
which he used to show that there is an abundance of singular
masas in separable ${\rm{II}}_1$ factors, \cite{Po_conv}.
Subsequently Popa, \cite{P3}, showed that a masa in a separable
${\rm{II}}_1$ factor is singular if, and only if, it has 
$\delta$--invariant $1$. An invariant $\alpha(A)$ for a masa $A$ in a
$\rm{II}_1$ factor, based on unitary perturbations of $A$, was
defined by the first two authors in \cite{SS1} with strong
singularity corresponding to $\alpha(A) =1$. In that paper the
masa $A$ was shown to be singular if $\alpha(A) > 0$.
Theorem~\ref{thm2.5} implies that a masa in a separable
${\rm{II}}_1$ factor is singular if, and only if, it is strongly
singular. This is the analogous result to Popa's
$\delta$--invariant one with unitaries in $N$ replacing his
nilpotent partial isometries, whose domains and ranges are
orthogonal projections in the masa. For a masa $A$ in a separable
$\rm{II}_1$ factor the results in \cite{P3} show that $\delta(A)$ is
either $0$ or $1$, and Theorem~\ref{thm2.5} implies that $\alpha(A)$ also takes
only these two values (see \cite{SS1} for the definition of $\alpha(A)$). 

The main result of the paper is Theorem~\ref{thm2.5}, the proof of
which is given in the lemmas that precede it. These in turn are
based on  results of Sorin Popa, \cite[Thm. 2.1, Cor. 2.3]{Po_ht}, which have
their origin in \cite{P4}. We also give two applications of our results. One
shows the singularity of tensor products of singular masas (Corollary
\ref{cor2.6}), while the other shows that singular masas can usually be
studied in the setting of separable algebras (Theorem \ref{thm10}).

Much of the work in this paper was accomplished at the Spring Institute on
Non-commutative Geometry and Operator Algebras, held May 9--20 2005 at
Vanderbilt University. The lecture series presented by Sorin Popa at this
conference provided the basis for our results below. It is a pleasure to
express our gratitude to the organizers and to the NSF for providing financial
support during the conference.

\newpage

\section{Main Results}\label{sec2}

Our first lemma is essentially contained in \cite[Corollary 2.3]{Po_ht}. The
proof will amount to identifying the subalgebras to which this
corollary is applied. 

\begin{lem}\label{lem2.3}
Let $A$  be a masa in a ${\rm{II}}_1$
factor $\N$ and let $e,f\in A$  be nonzero projections
with the property that no nonzero partial isometry $w\in \N$
satisfies the conditions
\begin{equation}\label{eq2.5a}
ww^*,ww^* \in A,\ ww^*\leq e, \ w^*w\leq f, {\text{ and }} w^*Aw=Aw^*w.
\end{equation}
 If $\vp >0$ and
$x_1,\ldots ,x_k \in \N$ are given,  then there exists a unitary $u\in A$
such that
\begin{equation}\label{eq2.6}
\|\bb E_{A}(fx_ieux_j^*f)\|_{2} <\vp
\end{equation}
for $1\leq i,j\leq k$.
\end{lem}

\begin{proof}
Define two subalgebras $B_0=Ae$ and $B=Af$ of $N$. The hypothesis implies the
negation of the fourth condition for $B_0$ and $B$ in \cite[Theorem 2.1]{Po_ht}
and so
\cite[Corollary 2.3]{Po_ht} can be applied. Thus, given $a_1,\ldots,a_k\in N$
and $\vp >0$, there exists a unitary $u_1\in B_0$ such that
\begin{equation}\label{eq.2.1}
\|{\mathbb E}_B(a_iu_1a_j^*)\|_{2}< \vp,\ \ \ 1\leq i,j\leq k.
\end{equation}
The result follows from this by taking $a_i=fx_ie$, $1\leq i \leq k$, replacing
$u_1$ by the unitary $u=u_1+(1-e)\in A$, and  replacing ${\mathbb E}_B$ in
\eqref{eq.2.1} by ${\mathbb E}_{A}$. Note that these two conditional
expectations
agree on $fNf$.
\end{proof}
Below, we will use the notation ${\mathcal U}(M)$ for the unitary group of any
von~Neumann algebra $M$. We will also need the well known fact that if $x\in M$
and $B$ is a von~Neumann subalgebra, then $\mathbb E_{B'\cap M}(x)$ is the
unique element of minimal norm in the $\|\cdot\|_2$--closure of
${\rm{conv}}\,\{uxu^*\colon u\in {\mathcal U}(B)\}$. We do not have an exact
reference for this, but it is implicit in \cite{Ch}.

\begin{lem}\label{lem2.4}
Let $A$ be a singular masa in a 
 ${\rm{II}}_1$ factor $\N$. If
$x_1,\ldots ,x_k\in \N$ and $\vp >0$ are
given, then there is a unitary $u\in
A$ such that
\begin{equation}\label{eq2.16}
\|\bb E_A(x_iux_j^*)-\bb E_A(x_i)u\bb E_A(x_j^*)\|_{2}<\vp
\end{equation}
for $1\leq i,j\leq k$.
\end{lem}
\begin{proof}
If $x,y\in \N$ and $u\in A$, then
\begin{equation}\label{eq2.17}
\bb E_A(xuy^*)-\bb E_A(x)u\bb E_A(y^*)=
\bb E_A((x-\bb E_A(x))u\bb E_A(y-\bb E_A(y))^*)
\end{equation}
by the module properties of $\bb E_A$. Thus (\ref{eq2.16})
follows if we can
establish that
\begin{equation}\label{eq2.18}
\|\bb E_A(x_iux_j^*)\|_{2}<\vp,\ \ \ 1\leq i,j\leq k,
\end{equation}
when the $x_i$'s also satisfy $\bb E_A(x_i)=0$
for $1\leq i \leq k$. We assume
this extra condition, and prove (\ref{eq2.18}). By scaling, there is no loss of
generality in assuming $\|x_i\|\leq 1$ for $1\leq i\leq k$.

Let $\delta =\vp /4$. In the separable case,  \cite{Po_conv} gives
 a finite dimensional
abelian
subalgebra $A_1 \subseteq A$ with minimal
projections $e_1,\ldots ,e_n$ such
that
\begin{equation}\label{eq2.19}
\|\bb E_{A_1'\cap\N}(x_i)-\bb E_A(x_i)\|_
{2}
 < \delta
\end{equation}
for $1\leq i \leq k$.
The assumption that $\bb E_A(x_i)=0$  allows us to rewrite \eqref{eq2.19} as
\begin{equation}\label{eq2.19a}
\|\bb E_{A_1'\cap\N}(x_i)\|_
{2}=
\|\bb E_A(x_i)-\bb E_{A_1'\cap\N}(x_i)\|_
{2} < \delta
\end{equation}
for $1\leq i \leq k$, leading to
\begin{equation}\label{eq2.19z}
\|\sum_{m=1}^n e_mx_ie_m\|_{2}=\|\bb E_{A_1'\cap\N}(x_i)\|_
{2}
 < \delta
\end{equation}
since $\sum_{m=1}^n e_mx_ie_m=\bb E_{A_1'\cap\N}(x_i)$.  Any partial isometry
$v\in \N$
satisfying $vAv^*=Avv^*$
has the form $pu$ for a projection $p\in A$ and a
normalizing unitary $u\in
\N$, \cite{Dye,JP}. The singularity of $A$ then
shows that $v\in A$, making it impossible to satisfy the
two
inequalities $vv^*\leq e_m$ and $v^*v\leq (1-e_m)$
simultaneously
unless $v=0$. Thus no nonzero partial isometry $v\in\N$
satisfies
$vv^*\leq e_m$, $v^*v\leq (1-e_m)$, and  $vAv^*=Avv^*$.
The hypothesis of
Lemma \ref{lem2.3} is satisfied, and applying
this result with $\vp$ replaced
by $\delta /n$ gives unitaries $u_m\in A$ such that
\begin{equation}\label{eq2.20}
\|\bb E_A((1-e_m)x_ie_mu_mx_j^*(1-e_m))\|_{2}<\delta /n
\end{equation}
for $1\leq m \leq n$ and $1\leq i,j\leq k$.
Define a unitary $u\in A$ by $u=\sum_{m=1}^n u_me_m$,
and let $y_i=\sum_{m=1}^n
(1-e_m)x_ie_m$, for $1\leq i\leq k$. We have
\begin{equation}\label{eq2.21}
x_i-y_i=x_i-\sum_{m=1}^n(1-e_m)x_ie_m
=\sum_{m=1}^n e_mx_ie_m
=\bb E_{A_1'\cap \N}(x_i)
\end{equation}
for $1\leq i\leq k$. The inequalities 
\begin{equation}\label{eq2.21a}
\|x_i-y_i\|_{2}
<\delta, 
\ \ \,\|x_i-y_i\|\leq \|x_i\|\leq 1, \ \ \, \|y_i\|\leq 2,
\end{equation}
for $1\leq i\leq k$,
follow immediately from \eqref{eq2.19a} and \eqref{eq2.21}.

 If we apply $\bb E_A$ to the identity
\begin{equation}\label{eq2.22}
x_iux_j^*=(x_i-y_i)ux_j^*+y_iu(x_j^*-y_j^*)+y_iuy_j^*,
\end{equation}
then \eqref{eq2.21a} 
gives
\begin{align}\label{eq2.23}
\|\bb E_A(x_iux_j^*)\|_{2}&\leq
\|x_i-y_i\|_{2}+\|y_i\|\|x_j-y_j\|_{2}
+ \|\bb E_A(y_iuy_j^*)\|_{2}\nonumber\\&<3\delta +
\|\bb E_A(y_iuy_j^*)\|_{2},
\end{align}
and we estimate the last term. The identity
\begin{equation}\label{eq2.24}
y_iuy_j^*=\sum_{m,s=1}^n(1-e_m)x_ie_mue_sx_j^*(1-e_s)
=\sum_{m=1}^n(1-e_m)x_ie_mu_mx_j^*(1-e_m)
\end{equation}
holds because each $e_s$ commutes with $u$ and $e_me_s=0$ for $m\ne s$.
The last sum has $n$ terms, so the
inequalities
\begin{equation}\label{eq2.25}
\|\bb E_A(y_iuy_j^*)\|_{2}<\delta,\ \ \ 1\leq i,j\leq k,
\end{equation}
are immediate from (\ref{eq2.20}). Together (\ref{eq2.23})
and (\ref{eq2.25})
 yield
\begin{equation}\label{eq2.26}
\|\bb E_A(x_iux_j^*)\|_{2}<3\delta +\delta =
\vp,\ \ \ 1\leq i,j\leq k,
\end{equation}
as required.

In the general case, we obtain $A_1$ and \eqref{eq2.19} as follows. Since $A$
is
a masa,
\begin{equation}
\mathbb E_{A'\cap N}(x_i)=\mathbb E_A(x_i)=0,\ \ \ 1\leq i\leq k.
\end{equation}
Now $\mathbb E_{A'\cap N}(x_i)$ is the element of minimal
$\|\cdot\|_2$--norm in the $\|\cdot\|_2$--closed convex hull of
$\{wx_iw^*\colon w\in {\mathcal U}(A)\}$, so we may select a finite number
of unitaries $w_1,\ldots , w_r\in A$ such that each set
$\Omega_i={\text{conv}}\,\{w_jx_iw_j^*\colon 1\leq j \leq r\}$, $1\leq
i\leq k$, contains an element whose $\|\cdot\|_2$--norm is less than
$\delta$. The spectral theorem allows us to make the further assumption
that each $w_j$ has finite spectrum, whereupon these unitaries generate a
finite dimensional subalgebra $A_1\subseteq A$. Then, for $1\leq i\leq k$,
$\mathbb E_{A_1'\cap N}(x_i)$ is the element of smallest norm in the
$\|\cdot\|_2$--closed convex hull of $\{wx_iw^*\colon w\in {\mathcal
U}(A_1)\}$, and since this set contains $\Omega_i$, we see that \eqref{eq2.19}
is valid in general.
\end{proof}

In \cite{SS1}, a masa $A$ in a  ${\rm{II}}_1$ factor $N$ was defined to
have the
{\emph{asymptotic homomorphism property}} ({\emph{AHP}}) if there exists a
unitary
$v\in A$ such that
\begin{equation}\label{eq5.1}
\lim_{|n|\to\infty}\|\bb E_A(xv^ny)-\bb E_A(x)v^n\bb E_A(y)\|_{2}=0
\end{equation}
for all $x,y\in N$. In that paper it was shown that strong singularity is a
consequence of this property. Subsequently it was observed
in \cite[Lemma 2.1]{RSS} that a
weaker property, which we will call the {\emph {weak asymptotic homomorphism
property}}, ({\emph{WAHP}}), suffices to imply strong singularity: given $\vp
>0$
and $x_1,\ldots ,x_k,y_1,\ldots ,y_k\in N$, there exists a unitary $u\in A$
such that
\begin{equation}\label{eq5.2}
\|\bb E_A(x_iuy_j)-\bb E_A(x_i)u\bb E_A(y_j)\|_{2}<\vp
\end{equation}
for $1\leq i,j\leq k$. Since the {\emph{WAHP}} is a consequence of applying
Lemma
\ref{lem2.4} to the set of elements
$x_1,\ldots ,x_k,y_1^*,\ldots ,y_k^*\in N$, we immediately obtain the main
result of the paper from these remarks:
\begin{thm}\label{thm2.5}
Let $A$ be a singular masa in a   ${\rm{II}}_1$ factor $N$. Then
$A$ has
the WAHP and is strongly singular.
\end{thm}

The following observation on the tensor product of masas may be known to some
experts, but we  have not found a reference.
\begin{cor}\label{cor2.6}
For $i=1,2$, let $A_i\subseteq N_i$ be masas in   ${\rm{II}}_1$
factors.
If $A_1$ and $A_2$ are both singular, then $A_1\overline\otimes A_2$ is also a
singular masa in $N_1\overline\otimes N_2$.
\end{cor}
\begin{proof}
Lemma \ref{lem2.4} and the remarks preceding Theorem \ref{thm2.5} show that
singularity and the {\emph{WAHP}} are equivalent for masas in  
${\rm{II}}_1
$ factors. By Tomita's commutant theorem, $A_1\overline\otimes A_2$ is a masa
in  $N_1\overline\otimes N_2$, and it is straightforward to verify that the
{\emph{WAHP}}
carries over to tensor products (see \cite[Proposition 1.4.27]{W}), since it
suffices to check this property on a 
$\|\cdot\|_
{2}$--norm dense set of elements, in this case the span of $\{x\otimes
y\colon x\in N_1,\,\,y\in N_2\}$.
\end{proof}

As an application of these results, we end by showing that the study of
singular masas can, in many instances, be reduced to the separable case.
The techniques of the proof have their origin in \cite[Section 7]{CPSS}.
\begin{thm}\label{thm10}
Let $N$ be a ${\rm{II}}_1$ factor with a singular masa $A$ and let $M_0$ be a
separable von~Neumann subalgebra of $N$. Then there exists a separable
subfactor $M$ such that $M_0\subseteq M \subseteq N$ and $M\cap A$ is a
singular masa in $M$.
\end{thm}

\begin{proof} We will construct $M$ as the weak closure of the union of an
increasing sequence $M_0\subseteq M_1\subseteq M_2\subseteq\ldots$ of
separable von~Neumann subalgebras, chosen by induction. These will have an
increasing sequence of abelian subalgebras $B_k\subseteq M_k$ with certain
properties.

For a von~Neumann algebra $Q\subseteq N$ and for any $x\in N$, $K_Q(x)$ will
denote the set ${\text{conv}}\,\{uxu^*\colon u\in {\mathcal U}(Q)\}$, and the
$\|\cdot\|$-- and $\|\cdot\|_2$--closures will be denoted $K_Q^n(x)$ and
$K_Q^w(x)$ respectively.  The inclusions $K_Q(x)\subseteq K_Q^n(x)\subseteq
K_Q^w(x)$ are immediate.  The induction hypothesis is: each $M_k$ is separable,
$M_k\subseteq M_{k+1}$, and for a fixed sequence 
$\{y_{k,r}\}_{r=1}^{\infty}$ in the unit ball of $M_k$ which is
$\|\cdot\|_2$--dense in the $\|\cdot\|_2$--closure of this ball,
\begin{itemize}
\item[(i)] $\mathbb E_A(y_{k,r})\in B_{k+1}\cap K_{B_{k+1}}^w(y_{k,r})$ for
$r\geq 1$, where $B_{k+1}=M_{k+1}\cap A$;

\item[(ii)] $K_{M_{k+1}}^n(y_{k,r})\cap {\mathbb C}1$ is nonempty for $r\geq
1$;

\item[(iii)] given $\vp >0$, $r\geq 1$ and a projection $p\in B_k$, there
exists $u\in {\mathcal U}(B_{k+1})$ such that
\[\|\mathbb E_A((1-p)y_{k,s}puy_{k,t}^*(1-p))\|_2< \vp \] for all $1\leq s,t
\leq r$. 
\end{itemize}

We first show that such a sequence of algebras leads to the desired conclusion.
Let $M$ and $B$ be respectively the weak closures of the unions of the
$M_k$'s and $B_k$'s. Since $K_{M_k}^n(x) \subseteq K_{M}^n(x)$ for all $x\in
M$ and $k\geq 1$, condition (ii) and a simple approximation argument show that
$K_{M}^w(x)$ contains a scalar operator for all $x$ in the unit ball of $M$,
and thus for all $x\in M$ by scaling. Since
$K_{M}^w(z)=\{z\}$ for any central element $z\in M$, this shows that $M$ is a
factor, separable by construction.

Now consider $x\in M$; scaling allows us to assume without loss of generality
that $\|x\|\leq 1$. Condition (i) shows that $\mathbb E_A(y_{k,r})\in
K_{B}^w(y_{k,r})$ for $k,r\geq 1$, and an approximation argument then shows
that  $\mathbb E_A(x)\in K_{B}^w(x)\subseteq K_{A}^w(x)$. Since $\mathbb
E_A(x)$ is the element of minimal $\|\cdot\|_2$--norm in $K_{A}^w(x)$, it also
has this property in $K_{B}^w(x)$. But this minimal element is 
$\mathbb E_{B'\cap M}(x)$, showing that $\mathbb E_{B'\cap M}(x)=\mathbb
E_{A}(x)$. If we further suppose that $x\in B'\cap M$, then $x=\mathbb
E_{A}(x)$.
Thus $B'\cap M\subseteq A$ and is abelian. Condition (i) also shows that $ 
\mathbb E_{A}(x)\in B$, and so $B'\cap M\subseteq B$. Since $B$ is abelian, we
have equality, proving that $B$ is a masa in $M$. We can now conclude that $
\mathbb E_{A}(y)=\mathbb E_{B}(y)$ for all $y\in M$, and that $B=M\cap A$.

Another $\|\cdot\|_2$--approximation argument, starting from (iii), gives
\begin{equation}
\inf\,\{\max_{1\leq i,j\leq r}\|\mathbb E_B((1-p)x_{i}pux_{j}^*(1-p))\|_2\colon
u\in {\mathcal U}(B)\}=0
\end{equation}
for an arbitrary finite set of elements $x_1,\ldots ,x_r\in M$, $\|x_i\|\leq
1$, and any
projection $p\in B$, noting that
$\mathbb E_A$ and $\mathbb E_B$ agree on $M$. By scaling, it is clear that this
equation holds generally without the restriction $\|x_i\|\leq 1$. We have now
established the conclusion
of Lemma \ref{lem2.3} for $B$ from which singularity of $B$ follows, as in the
proof of Lemma \ref{lem2.4}. It remains to construct the appropriate
subalgebras $M_k$.

To begin the induction, let $B_0=M_0\cap A$, and suppose that $B_k\subseteq
M_k$ have been constructed. Consider a fixed  sequence
$\{y_{k,r}\}_{r=1}^{\infty}$ in the unit ball of $M_k$, $\|\cdot\|_2$--dense in
the $\|\cdot\|_2$--closure of this ball. The Dixmier approximation theorem,
\cite[Theorem 8.3.5]{KR}, allows
us to obtain a countable number of unitaries, generating a separable
subalgebra $Q_0\subseteq N$, so that $K_{Q_0}^n(y_{k,r})\cap {\mathbb C}1$ is
nonempty for $r\geq 1$. Let $\{p_m\}_{m=1}^{\infty}$ be a sequence which is
$\|\cdot\|_2$--dense in the set of all projections in $B_k$. The singularity of
$A$ ensures that the hypothesis of Lemma \ref{lem2.3} is met when $e=p_m$ and
$f=1-p_m$. Thus, for integers $m,r,s\geq 1$, there is a unitary $u_{m,r,s}\in
A$ such that
\begin{equation}
 \|\mathbb E_A((1-p_m)y_{k,i}p_mu_{m,r,s}y_{k,j}^*(1-p_m))\|_2<1/s
\end{equation}
for $1\leq i,j\leq r$. These unitaries generate a separable von~Neumann algebra
$Q_1\subseteq A$, and an approximation argument establishes (iii) provided that
$Q_1 \subseteq B_{k+1}$.

Since $\mathbb E_A(y_{k,r})$ is the minimal $\|\cdot\|_2$--norm element in
$K_A^w(y_{k,r})$, we may find a countable number of unitaries generating a
separable subalgebra $Q_2\subseteq A$ so that $\mathbb E_A(y_{k,r})\in 
K_{Q_{2}}^w(y_{k,r})$ for $r\geq 1$. The proof is completed by letting
$M_{k+1}$ be the separable von~Neumann algebra generated by $M_k$, $\mathbb
E_A(M_k)$, $Q_0$, $Q_1$ and $Q_2$. The subspaces $\mathbb E_A(M_k)$ and $Q_2$  
are included to ensure that (i) is satisfied, (ii) holds by the choice of
$Q_0$, and $Q_1$ guarantees the validity of (iii).
\end{proof}

\end{document}